\documentclass[a4paper,10pt]{article}
\usepackage{amssymb,amsmath,amsthm,latexsym}
\usepackage{color}
\usepackage{textcomp}
\usepackage[active]{srcltx}
\usepackage{colortbl}

\usepackage{hyperref}

\usepackage{color}
\usepackage{mathrsfs}
\usepackage{graphicx}
\usepackage{lineno}
\usepackage{adjustbox}
\usepackage{rotating}

\usepackage{longtable}

\usepackage{tikz}
\usetikzlibrary{shapes.geometric, arrows}

\tikzstyle{process} = [rectangle, minimum width=2.8cm, minimum height=0.8cm, text centered, draw=black, font=\small]
\tikzstyle{decision} = [ellipse, minimum width=2.5cm, minimum height=0.8cm, text centered, draw=black, font=\small]
\tikzstyle{arrow} = [thick,->,>=stealth]

\usepackage[multiple]{footmisc}

\usepackage{fancyhdr}

\newcommand{\keywords}[1]{%
  \small
  \textbf{\textit{Keywords---}} #1
}

\newcommand{\mscclass}[1]{\par\noindent\textbf{MSC:} #1}

\newtheorem{thm}{Theorem}[section]

\newtheorem{prop}[thm]{Proposition}
\theoremstyle{definition}
\newtheorem{defn}[thm]{Definition}
\theoremstyle{remark}

\numberwithin{equation}{section}
\newtheorem{exmp}[thm]{Example}

\title{Salem numbers less than $49/37$}

\author{J.-M. Sac-\'Ep\'ee\footnote{J.-M. Sac-\'Ep\'ee, IECL, Université de Lorraine,
 France, jean-marc.sac-epee@univ-lorraine.fr}}

\begin{document}

\maketitle
\begin{center}
\emph{Dedicated to the memory of Professor Georges Rhin}
\end{center}
\newcommand{\D}{\mathbb{D}}
\newcommand{\C}{\mathbb{C}}
\newcommand{\N}{\mathbb{N}}
\newcommand{\R}{\mathbb{R}}
\newcommand{\Z}{\mathbb{Z}}
\newcommand{\dist}{\operatorname{dist}}

\renewcommand{\qedsymbol}{$\blacksquare$}

\makeatletter
\def\blfootnote{\xdef\@thefnmark{}\@footnotetext}
\makeatother

\begin{abstract}
A certain number of lists of \emph{small} Salem numbers, some of which are certified as complete, are available online. Notably, the website of M. J. Mossinghoff features a list of 47 Salem numbers smaller than 1.3 (\cite{MossinghoffList}), as well as complete lists of Salem numbers of fixed degrees that are smaller than various bounds (\cite{MossinghoffList2}). The objective of this work is to advance the understanding of Salem numbers by providing a list of Salem numbers smaller than a threshold of $49/37$ ($\approx 1.324324$) through the implementation of a method based on integer linear programming and uniform sampling of root separators. Beyond the intrinsic interest of the newly detected Salem numbers, the rediscovery of already known Salem numbers through this alternative method could offer valuable insights into the potential completeness of already existing lists for degrees where it has not yet been proven.
\end{abstract}
\mscclass{11K16}


\section{Introduction} 
Salem numbers have drawn significant interest since at least the 1940s. The focus on \emph{small Salem numbers} has been particularly prominent, partly due to their connection with Lehmer's problem, which concerns integer polynomials with small Mahler measure. Numerous researchers have investigated these numbers, with efforts tracing back at least to D. W. Boyd's work in the 1970s. In the following lines, we provide a brief overview of some key properties of Salem numbers.
\subsection{Definitions and examples}
\begin{defn}
Let $\alpha \in \mathbb{C}$.
 \begin{itemize}
    \item  We say that $\alpha$ is an \emph{algebraic number} if $\alpha$ is a root of a non-zero polynomial in $\mathbb{Z}[x]$, the ring of polynomials in one indeterminate $x$ with coefficients in $\mathbb{Z}$.
    \item We say that $\alpha$ is an \emph{algebraic integer} if $\alpha$ is a root of a \emph{monic} polynomial in $\mathbb{Z}[x]$.
    \item The \emph{degree} of $\alpha$ is defined as the degree of its minimal polynomial.
 \end{itemize}
\end{defn}
\begin{exmp}
The algebraic integer $i$ has degree 2, because its minimal polynomial is $x^2 + 1$.
\end{exmp}

\begin{defn}
A \emph{Salem number} is an algebraic integer $\tau$ such that:
\begin{itemize}
    \item $\tau$ is a real number strictly greater than $1$,
    \item All of its Galois conjugates (i.e., the other roots of its minimal polynomial) have absolute value less than or equal to $1$,
    \item At least one of its Galois conjugates has absolute value equal to $1$.
\end{itemize}
\end{defn}
\begin{defn}
The minimal polynomial of a Salem number is called a \emph{Salem polynomial}.
\end{defn}
\begin{exmp}
The Salem number $\tau_0 = 1.176280\dots$ has the minimal polynomial $P_0(x) = x^{10} + x^9 - x^7 - x^6 - x^5 - x^4 - x^3 + x + 1$, known as the \emph{Lehmer polynomial}. The real number $\tau_0$ is the smallest known Salem number to date (\cite{Lehmer1933}).
\end{exmp}

\subsection{Some properties}

\begin{prop}
\
\begin{itemize}
    \item The minimal polynomial of a Salem number $\tau$ is reciprocal and has even degree greater than or equal to 4.
    \item The real number $\displaystyle\frac{1}{\tau}$ is one of the Galois conjugates of $\tau$.
    \item All the roots of the minimal polynomial of $\tau$, except for $\tau$ and $\displaystyle\frac{1}{\tau}$, have absolute value equal to $1$.
\end{itemize}
\end{prop}

A very rich and comprehensive reference on Salem numbers is \cite{Smyth2015}.

\subsection{Small Salem numbers}
Salem numbers less than $1.3$ are referred to as \emph{small} in the literature. In 1977, D. W. Boyd provided a list of $39$ such numbers in \cite{Boyd1977}. The following year, D. W. Boyd added $4$ more small Salem numbers in \cite{Boyd1978}, before $4$ additional small Salem numbers were discovered and published by M. J. Mossinghoff in 1993 (\cite{Mossinghoff1993}). Notably, among the $4$ values provided by M. J. Mossinghoff, one has a degree of $46$, which is the highest known degree for a small Salem number.

The compilation of these sources thus yields a list of $47$ small Salem numbers, which is available on M. J. Mossinghoff's website (\cite{MossinghoffList}). This list provides, for each Salem number, its degree and the coefficients of its minimal polynomial. Naturally, only the relevant part of the coefficients is given, as these polynomials are reciprocal.

This list was first certified complete up to degree $40$ in \cite{FlammangGrandcolasRhin1999}, and later up to degree $44$ in \cite{MossinghoffRhinWu2008}. Only the degree $46$ polynomial mentioned earlier does not belong to the certified complete list.

In this work, we propose to extend the list \cite{MossinghoffList} by providing Salem numbers that are smaller than $\eta = 49/37 $. The value chosen for $\eta$ might seem somewhat arbitrary and therefore deserves a brief justification. Let us point out that the real number $\eta$ is close to the plastic constant, the real root of $x^3-x-1$, whose approximate value is $1.324718$. The plastic constant holds particular significance in the study of Salem numbers. Let us recall that a Pisot-Vijayaraghavan number is a real algebraic integer greater than $1$, all of whose Galois conjugates are less than $1$ in absolute value. R. Salem provided a construction which, given a Pisot-Vijayaraghavan number, produces two sequences of Salem numbers, one approaching the Pisot-Vijayaraghavan number from below and the other from above (\cite{Salem1963}, \cite{Smyth2015}). In turn, D. W. Boyd proved that every Salem number appears as a member of one of the sequences of Salem numbers arising from Pisot-Vijayaraghavan numbers (see \cite{Boyd1977} or \cite{Smyth2015}). It is worth mentioning that the plastic constant, which is the smallest Pisot-Vijayaraghavan number, is the smallest known limit point of the Salem numbers. From the previous lines, we deduce that there is an infinite number of Salem numbers smaller than the plastic constant, and it is possible to explicitly construct a sequence of such numbers. This explains why the value $1.3$ has traditionally been chosen as the threshold for defining \emph{small} Salem numbers.  The choice in this work is to set a new threshold exceeding $1.3$, while remaining strictly below the plastic constant, for clear reasons mentioned above. Our choice fell on the value $49/37$, which has an approximate value of $1.324324$.

To construct our extended list, we used a variation of a method based on integer linear optimization, which has allowed us to obtain interesting results in various contexts. Notably, we followed this approach to provide, in \cite{ElOtmaniMaulRhinSac-Epee2014}, a list of degree $16$ irreducible monic polynomials with integer coefficients, all real and strictly positive roots and minimal trace (specifically, with trace $29$). Similarly, we provided examples of Salem numbers of degree $34$ and trace $-3$ in \cite{ElOtmaniRhinSac-Epee2015} and \cite{Sac-Epee2024}.

Here, it is not the trace of Salem numbers that interests us, but their size. Our approach, based on integer linear optimization, is well-suited for addressing this problem. In the following lines, we only outline the method, as it is described in detail in the aforementioned articles.

\section{Algorithmic approach}
\subsection{Reduction to a half-degree problem}

The minimal polynomial $P$ of a Salem number $\tau$, as mentioned above, is a reciprocal polynomial of even degree (say $2d$). It can be written as
\[
P(x) = x^{2d} + c_1 x^{2d-1} + c_2 x^{2d-2} + \ldots + c_{d-1} x^{d+1} + c_d x^d + c_{d-1} x^{d-1} + \ldots + c_2 x^{2} + c_1 x + 1,
\]
and can be rewritten as
\[
P(x) = x^d Q\left(x + \frac{1}{x}\right),
\]
where $Q$ is a degree $d$ monic polynomial of $\mathbb{Z}[x]$, whose roots are all real numbers and lie within the interval $(-2,2)$, except for one root that is strictly greater than $2$. The polynomial $Q$ is obtained explicitly as follows: the first coefficients of $P$ are the integers $1, c_1, c_2, \ldots, c_d$ (the other coefficients of $P$ are determined by symmetry), and then
\[
Q(x) = C_{d+1}(x) + c_1 C_d(x) + c_2 C_{d-1}(x) + \ldots + c_d C_1(x),
\]
where the $C_i$'s are the Chebyshev polynomials associated with the interval $(-2,2)$. As a reminder, we have

\begin{itemize}
  \item[] $C_1(x) = 1$
  \item[] $C_2(x) = x$
  \item[] $C_3(x) = x^2-2$
  \item[] $C_4(x) = x^3-3x$
  \item[] $C_5(x) = x^4-4x^2+2$
  \item[] $C_6(x) = x^5-5x^3+5x$
   \item[] $\vdots$
\end{itemize}
For example, the Lehmer polynomial $P_0(x) = x^{10} + x^9 - x^7 - x^6 - x^5 - x^4 - x^3 + x + 1$ can be written as
\[
P_0(x) = C_{6}(x) + c_1 C_5(x) + c_2 C_{4}(x) + c_3 C_{3}(x) + c_4 C_{2}(x) + c_5 C_{1}(x) = x^5 Q_0\left(x + \frac{1}{x}\right),
\]
where $Q_0(x) = x^5 + x^4 - 5 x^3 - 5 x^2 + 4 x + 3$. It is easy to verify that $Q_0$ has $4$ roots in $(-2, 2)$ and a single root strictly greater than $2$.

Since $P$ can be easily obtained from $Q$ and vice versa, the search for Salem numbers of degree $2d$, which is equivalent to the search for Salem polynomials of degree $2d$, can be reduced to the search for monic polynomials of degree $d$, with integer coefficients, and all real roots, where $d-1$ roots lie in $(-2, 2)$ and the last one is strictly greater than $2$. In this work, we have an additional constraint: we want the polynomials $P$ (of degree $2d$) to be Salem polynomials, and furthermore, we want their unique real root greater than $1$ to also be strictly less than $\eta$. When transitioning from the search for $P$ to the search for $Q$, this constraint translates to the fact that the unique root of $Q$ greater than $2$ must also be strictly less than $\eta + \displaystyle\frac{1}{\eta}$.

\subsection{Reformulation as an Optimization Problem}
As mentioned earlier, we have an efficient algorithmic approach based on integer linear optimization for the search for polynomials of type $Q$, which we outline below.

Consider an irreducible polynomial in $\mathbb{Z}[x]$, $q(x) = x^d + \displaystyle\sum_{i=0}^{d-1} a_i x^i$, with all real roots, where $d-1$ roots lie in $(-2, 2)$ and a single root is strictly greater than $2$ and less than $\eta + \displaystyle\frac{1}{\eta}$. For such a polynomial $q$, any $(d+1)$-tuple $\left(-2, \beta_1, \dots, \beta_{d-2}, 2, \eta + \displaystyle\frac{1}{\eta}\right)$ such that $\beta_1, \dots, \beta_{d-2}$ and $2$ separate its $d$ roots, satisfies the condition that $q(-2)$, $q(\beta_1)$, $q(\beta_2)$, $\dots$, $q(\beta_{d-2})$, $q(2)$, $q\left(\eta + \displaystyle\frac{1}{\eta}\right)$ alternate in sign, which provides $d+1$ \emph{linear constraints} with respect to the \emph{integer} variables $a_i$. These linear constraints form an integer linear optimization problem whose feasibility we must investigate. In other words, we must test the existence of integers $a_i$ that satisfy all the linear constraints. The idea is to randomly and independently sample a large number of $(d-2)$-tuples $(\beta_1, \dots, \beta_{d-2})$ within the interval $(-2, 2)$ and use the $(d+1)$-tuples $\left(-2, \beta_1, \dots, \beta_{d-2}, 2, \eta + \frac{1}{\eta}\right)$ as root separators for the potential relevant polynomials. Solving the integer linear optimization problem (which we do for each sampled $(\beta_1, \dots, \beta_{d-2})$) does not involve maximization or minimization, as there is no objective function to optimize. The algorithm is reduced to a feasibility study for each randomly drawn $(d-2)$-tuple.

\subsection{Practical implementation}
In practice, the $\beta_i$ points are drawn uniformly from $(-2,2)$ in each trial. Once a $(d-2)$-tuple $(\beta_1, \dots, \beta_{d-2})$ is drawn, we solve the problem. If no $d$-tuple $(a_0, \dots, a_{d-1})$ is found that satisfies all the constraints, we immediately draw another sample. Conversely, if such a $d$-tuple $(a_0, \dots, a_{d-1})$ is found, this means we have found a polynomial $q$ (of degree $d$) of the type $Q$, which allows us to reconstruct a Salem polynomial $P$ of the desired type. If $q$ is irreducible in $\mathbb{Z}[x]$, we add the polynomial $P$ reconstructed from $q$ to our list. If $q$ is not irreducible, we keep the factor of $q$, say $q_1$, that has a root between $2$ and $\eta + \displaystyle\frac{1}{\eta}$, and we reconstruct a Salem polynomial from $q_1$. The following flowchart summarizes the approach used:
\newline

\begin{tikzpicture}[node distance=1.3cm]  

\node (start) [process] {Sampling of the $\beta_i$};
\node (test) [process, below of=start, yshift=-0.5cm] {Search for the $a_i$};
\node (ai_notfound) [decision, left of=test, xshift=-2.5cm] {$a_i$ not found};  
\node (ai_found) [decision, right of=test, xshift=3.5cm] {$a_i$  found};  
\node (exam) [process, right of=start, xshift=3.5cm] {Irreducibility of $q$};  

\draw [arrow] (start) -- (test);
\draw [arrow] (test) -- (ai_notfound);
\draw [arrow] (test) -- (ai_found);  
\draw [arrow] (ai_found) -- (exam);  
\draw [arrow] (exam) -- (start);  
\draw [arrow] (ai_notfound) |- (start);  

\end{tikzpicture}

We performed calculations related to the search for degree $d$ polynomials with real roots for all values of $d$ between $10$ and $32$. From the perspective of Salem polynomials, this means we searched for such polynomials up to degree $2d=64$.

For small values of $d$, we carried out $800,000$ uniform random draws (this number was chosen because we had access to an $80$-core server), and we progressively increased the number of draws to $1,600,000$ for larger values of $d$.

\section{Technical details and results}
\subsection{Technical details}

The program was entirely written in the \emph{Julia} language. Polynomial manipulations were performed using the \emph{Polynomials} library, while irreducibility tests were carried out with the \emph{Nemo} library.

The integer linear optimization problems for each sample were solved using the \emph{Gurobi} library (\cite{Gurobi}).

The results were obtained on a multiprocessor server featuring two Intel Xeon Gold 6138 CPUs, each with 20 cores (2 threads per core), totaling 80 logical processors (hyperthreading is enabled).

\subsection{Results}

The implementation of the algorithmic approach presented in this work found all known Salem numbers less than $1.3$, but no new numbers less than this value. In the following table, we therefore only present the Salem numbers between $1.3$ and $\eta$ resulting from our calculations. To make the table easy to read, we follow the presentation chosen in \cite{MossinghoffList}, in which the different columns respectively show the degrees of the Salem numbers, the numbers themselves, and the relevant coefficients of the minimal polynomials (since these are reciprocal). Up to degree $24$, the Salem numbers found using our method are already known and are listed on M. J. Mossinghoff's website (\cite{MossinghoffList}). These polynomials appeared in a $2014$ article by K. G. Hare and M. J. Mossinghoff (\cite{HareMossinghoff2014} ), in a study of specific representations of the sequences of Salem numbers associated with small Pisot-Vijayaraghavan numbers. As for us, we provide 10 new Salem numbers less than $49/37$, with degrees ranging from $26$ to $44$. They are highlighted in bold in the table above.
\begin{longtable}{|c|c|l|}
\hline
\textbf{2d} & \textbf{Salem numbers} & \textbf{Coefficients of Minimal Polynomial} \\
\hline
\endfirsthead

\hline
\textbf{2d} & \textbf{Salem numbers} & \textbf{Coefficients of Minimal Polynomial} \\
\hline
\endhead

\hline
\endfoot

12 & 1.302268805094& 1 -1 0 0 0 -1 1 \\
\textbf{32} & \textbf{1.302721444014} & \textbf{1 -1 0 -1 0 1 0 0 0 -1 1 -1 1 0 0 0 -1} \\
\textbf{32} & \textbf{1.303283348964} & \textbf{1 1 0 -1 -2 -2 -1 0 1 1 0 -1 -1 -1 0 1 1} \\
\textbf{30} & \textbf{1.303385419369} & \textbf{1 -1 0 0 -1 0 0 0 1 0 0 1 -1 0 0 -1} \\
\textbf{26} & \textbf{1.304697625411} & \textbf{1 0 -1 -1 0 0 0 0 0 0 0 1 0 -1} \\
22 & 1.305131378642& 1 -1 0 0 -1 0 0 1 0 0 0 -1 \\
\textbf{38} & \textbf{1.306473537533} & \textbf{1 -2 1 0 0 0 0 0 -1 1 0 -1 1 0 -1 0 1 0 -1 1} \\
\textbf{44} & \textbf{1.308071085577} & \textbf{1 0 1 -1 0 -2 -1 -2 -1 -1 0 0 1 0 1 0 1 0 1 0 1 0 1} \\
16 & 1.308409006213& 1 1 0 -1 -1 -1 -1 -1 -1 \\
\textbf{36} & \textbf{1.308966300288} & \textbf{1 0 1 -1 0 -2 -1 -2 -1 -1 0 0 1 0 1 0 1 0 1}\\
22 & 1.310180863375& 1 0 0 -1 -1 -1 -1 0 0 1 1 1 \\
22 & 1.312566632631& 1 -1 -1 1 0 0 0 -1 0 1 0 -1 \\
16 & 1.312773239526& 1 0 0 -1 0 -1 0 -1 0 \\
12 & 1.315914431926& 1 0 0 0 -1 -1 -1 \\
\textbf{40} & \textbf{1.316069252718} & \textbf{1 -1 0 0 -1 0 1 -1 1 0 -1 0 0 -1 1 0 0 1 0 -1 1} \\
14 & 1.318197504432& 1 -1 0 -1 1 0 0 -1 \\
20 & 1.319869661883& 1 1 0 -1 -1 -1 -1 -1 -1 -1 -1 \\
14 & 1.321101848259& 1 -1 0 0 -1 0 1 -1 \\
22 & 1.322014239618& 1 1 0 -1 -1 -1 -1 -1 -1 -1 -1 -1 \\
22 & 1.322692457903& 1 0 0 -1 0 -1 0 -1 0 -1 0 -1 \\
18 & 1.323198173512& 1 -1 -1 1 0 0 0 -1 0 1 \\
18 & 1.323576201647& 1 0 0 0 0 -1 -1 -1 -1 -2 \\
\textbf{26} & \textbf{1.323859346186} & \textbf{1 1 0 -1 -1 -1 -1 -1 -1 -1 -1 -1 -1 -1} \\
24 & 1.324071761641& 1 -1 0 0 0 -1 1 -1 0 0 0 -1 1 \\
\textbf{28} & \textbf{1.324231319862}& \textbf{1 1 0 -1 -1 -1 -1 -1 -1 -1 -1 -1 -1 -1 -1} \\
\end{longtable}

The fact that all known Salem numbers less than 1.3, as well as all Salem numbers less than $49/37$ with a degree smaller than 24, have been recovered using our method seems highly significant to us. This is not only a testament to the effectiveness of our approach but also provides a basis for assessing the completeness of the lists available online. Indeed, two factors must be considered. The first is that our method is based on uniform random sampling, meaning it is free from any bias in its ability to detect values. The second is that not only were the already known values rediscovered using this method, but they were also found repeatedly, many times over. Similarly, the newly discovered values appeared numerous times throughout the program’s execution. Thus, while this is not theoretically proven, there are reasonable grounds to believe that for the degrees examined, we have likely identified all that could be found. The same argument leads us to think that the list of Salem numbers less than 1.3, certified complete up to degree 44, is likely complete up to significantly higher degrees, given that all its values were repeatedly detected without any new ones emerging.

\section{Conclusion}
A new Salem number smaller than $1.3$ would of course have been welcome. In any case, our line of research has proven fruitful, as $10$ new Salem numbers lower than $49/37$  are now available. It should be noted, of course, that this algorithmic approach does not  guarantee the exhaustiveness of the table it provides. But since it is based on random sampling, the repeated discovery of values already present in our list suggests that any missing values are likely to be few. In particular, the repeated discovery of known Salem numbers below $1.3$ using this algorithmic approach, without finding new ones, suggests that List \cite{MossinghoffList} might be complete well beyond degree $44$ certified in \cite{MossinghoffRhinWu2008}.

\section{Acknowledgement}
The author of this article wishes to express deep gratitude to the anonymous referee, whose valuable advice has greatly improved the final manuscript. They are sincerely thanked for it.


\begin{thebibliography}{<n>}

\bibitem{Smyth1971} C. J. Smyth, On the product of the conjugates outside the unit circle of an algebraic integer, Bulletin of the London Mathematical Society, vol. 3, no. 2, pp. 169-175, 1971.
\bibitem{MossinghoffList2} M. J. Mossinghoff, \href{http://wayback.cecm.sfu.ca/~mjm/Lehmer/lists}{http://wayback.cecm.sfu.ca/~mjm/Lehmer/lists}
\bibitem{MossinghoffList} M. J. Mossinghoff, \href{http://wayback.cecm.sfu.ca/~mjm/Lehmer/lists/SalemList.html}{http://wayback.cecm.sfu.ca/~mjm/Lehmer/lists/SalemList.html}
\bibitem{MossinghoffRhinWu2008} M. J. Mossinghoff, G. Rhin, Q. Wu, Minimal Mahler measures, Experiment. Math. 17 (2008), no. 4, 451-458.
\bibitem{Smyth2015} C. J. Smyth, Seventy years of Salem numbers, C. Bull. Lond. Math. Soc. 47, No. 3, 379-395 (2015).
\bibitem{Boyd1977} D. W. Boyd, Small Salem numbers. Duke Math. J. 44 (1977), no. 2, 315-328.
\bibitem{Boyd1978} D. W. Boyd, Pisot and Salem numbers in intervals of the real line. Math. Comp. 32 (1978), no. 144, 1244-1260.
\bibitem{Mossinghoff1993} M. J. Mossinghoff, Polynomials with small Mahler measure. Math. Comp. 67 (1998), no. 224, 1697-1705, S11-S14.
\bibitem{FlammangGrandcolasRhin1999} V. Flammang, M. Grandcolas, G. Rhin, Small Salem numbers. Number theory in progress, Vol. 1 (Zakopane-Ko scielisko, 1997), 165-168, de Gruyter, Berlin, 1999.
\bibitem{MossinghoffRhinWu2008} M. J. Mossinghoff, G. Rhin, Q. Wu, Minimal Mahler measures. Experiment. Math. 17 (2008), no. 4, 451-458.
\bibitem{ElOtmaniMaulRhinSac-Epee2014} S. El Otmani, A. Maul, G. Rhin, J.-M. Sac-Epée, Finding degree-16 monic irreducible integer polynomials of minimal trace by optimization methods, Exp. Math. 23 (2014), no. 1, 1-5.
\bibitem{ElOtmaniRhinSac-Epee2015} S. El Otmani, G. Rhin, J.-M. Sac-Épée, A Salem number with degree 34 and trace $-3$, J. Number Theory 150 (2015), 21-25.
\bibitem{Sac-Epee2024} J.-M. Sac-Épée, Small degree Salem numbers with trace -3, Control Cybern. 52, No. 3, 335-346 (2023). 
\bibitem{Gurobi} Gurobi Optimization, LLC,  Gurobi optimizer reference manual, 2024, \newline \href{http://www.gurobi.com}{http://www.gurobi.com}
\bibitem{Salem1963} R. Salem, Algebraic numbers and Fourier analysis, D. C. Heath and Company, Boston, (1963).
\bibitem{Boyd1977} D. W. Boyd, Small Salem numbers, Duke Math. J. 44 (1977), no. 2, 315-328.
\bibitem{HareMossinghoff2014} K. G. Hare, M. J. Mossinghoff, Negative Pisot and Salem numbers as roots of Newman polynomials, Rocky Mountain J. Math. 44 (2014), no. 1, 113-138.
\bibitem{Lehmer1933} D. H. Lehmer, Factorization of certain cyclotomic functions, Ann. of Math. (2) 34 (1933), no. 3, 461-479.
\end{thebibliography}
\end{document}